\documentclass[10pt]{article}

\usepackage{theorem,amssymb,amsmath}

\advance \topmargin by -\headheight
\advance \topmargin by -\headsep
\evensidemargin \oddsidemargin
\author{J.-P. Allouche \\
CNRS, Inst.\ Math.\ Jussieu-PRG, UPMC \\
Case 247, 4 Place Jussieu \\
F-75252 Paris Cedex 05 France \\
{\tt jean-paul.allouche@imj-prg.fr}
}

\title{Two exercises of Comtet and two identities of Ruehr}

\date{ }

\def \proof{\bigbreak\noindent{\it Proof.\ \ }}

\def \endpf{{\ \ $\Box$ \medbreak}}

\newtheorem{theorem}{Theorem}

\newtheorem{corollary}{Corollary}
\newtheorem{proposition}{Proposition}
\theorembodyfont{\rm}

\newtheorem{remark}{Remark}

\begin{document}

\maketitle

\begin{abstract}
A question proposed by Kimura and proved by Ruehr, Kimura and others in 1980 states that
for any function $f$ continuous on $[-\frac{1}{2}, \frac{3}{2}]$ one has
$$
\int_{-1/2}^{3/2} f(3x^2 - 2x^3) dx = 2 \int_0^1 f(3x^2 - 2x^3) dx.
$$
In his proof Ruehr indicates, without giving an explicit proof, that this identity, applied to 
$f(t) = t^n$, implies two identities involving binomial sums, namely (after correction of a misprint)
$$
\sum_{0 \leq j \leq n} 3^j {3n-j \choose 2n} = 
\sum_{0 \leq j \leq 2n} (-3)^j {3n-j \choose n}
\ \
\mbox{\rm and}
\ \
\sum_{0 \leq j \leq n} 2^j {3n+1 \choose n-j} =
\sum_{0 \leq j \leq 2n} (-4)^j {3n+1 \choose n+1+j}.
$$
Using two identities given in a book of Comtet we provide an easy explicit way of deducing 
these identities from the above equality between integrals. Our derivation shows a link with the 
incomplete beta function, the binomial distribution law, the negative binomial distribution law, 
and a lemma used in a proof of a very weak form of the $(3x+1)$-conjecture.
\end{abstract}

\section{The Kimura-Ruehr identity and its combinatorial consequences}

In the section ``Problems and Solutions'' of a 1980 issue of the Monthly \cite{kimura-ruehr}, 
Kimura proposed and Ruehr, Kimura and others proved, a curious result involving definite
integrals, namely: if  $f$ is continuous on $[-\frac{1}{2}, \frac{3}{2}]$, then
$$
\int_{-1/2}^{3/2} f(3x^2 - 2x^3) dx = 2 \int_0^1 f(3x^2 - 2x^3) dx.
$$
In his (trigonometric) proof Ruehr writes that this equality is equivalent to all such equalities where 
$f(t) = t^n$, and $n$ any nonnegative integer. He then writes ``The next natural step, introducing
binomial expansions and integrating term by term, leads, however to several combinatorial
identities such as ~(...) which are apparently new and no easier to prove than the original
equations. Their validity, of course, follows from the trigonometric proof''. The combinatorial
identities are indicated in the abstract above (after correction as indicated in \cite{mtwz},
also see \cite{alzer-prodinger}).The author of the present paper tried to deduce the combinatorial
identities from the identity between definite integrals, then he asked colleagues: no immediate
proof was found. Other papers on these combinatorial identities give some generalizations, but 
none gives the ``missing'' proof. The purpose of the present paper is to provide an explicit 
reasonably simple proof that puts to light two exercises in the famous book of Comtet on 
Combinatorics. Furthermore it shows a relationship with the incomplete beta function, the binomial 
distribution law, the negative binomial distribution law, and a lemma used in a proof of... a very 
weak form of the $(3x+1)$-conjecture.

\section{Two exercises from Comtet's book on combinatorics}

The nice books of Comtet on Combinatorics exist in two versions, a French one in two volumes
and an English one. Among the multitude of results in these books, there are in the French version 
\cite{comtet-fr} two exercices (Exercise 12 p.~91 and Exercise 35 p.~179) giving (classical?) 
equalities relating certain binomial sums and integrals. The first exercise is present in the English 
version \cite{comtet-eng} (Exercise~12 p.~76), while the second exercise seems to be missing 
(except for the particular case $x=1/2$, see bottom of Page~72). We state these two exercices in 
the following proposition.

\begin{proposition}\label{exercises-comtet}
\begin{equation}\label{comtet1}
\sum_{0 \leq i \leq k} {n \choose i} a^{n-i} b^i =
(n-k){n \choose k} \int_b^{a+b} t^k (a+b-t)^{n-k-1} dt
\end{equation}
and
\begin{equation}\label{comtet2}
\sum_{m \leq k \leq n} {k-1 \choose m-1} x^m (1-x)^{k-m} =
\sum_{m \leq k \leq n} {n \choose k} x^k (1-x)^{n-k}.
\end{equation}
\end{proposition}

\proof
To keep this paper self-contained we offer elementary proofs for these two equalities.
\begin{itemize}

\item{ } Equality~(\ref{comtet1}) is proved by induction on $k$. It is clearly true for $k=0$.
Suppose it is true for $k$, then, integrating by parts, one has
\begin{multline*}
(n-(k+1)){n \choose k+1} \int_b^{a+b} t^{k+1} (a+b-t)^{n-(k+1)-1} dt \\
\begin{array}{lll}
&=&
\displaystyle{n \choose k+1} \left( b^{k+1} a^{n-k-1} + \int_b^{a+b} (k+1) t^k (a+b-t)^{n-k-1} dt \right) \\
&=& 
\displaystyle{n \choose k+1} b^{k+1} a^{n-k-1} + (n-k) {n \choose k} \int_b^{a+b} t^k (a+b-t)^{n-k-1} \\
&=&
\displaystyle{n \choose k+1} b^{k+1} a^{n-k-1} + \sum_{0 \leq i \leq k} {n \choose i} a^{n-i} b^i 
\ \ \mbox{\rm (from the induction hypothesis)} \\
&=& 
\displaystyle\sum_{0 \leq i \leq k+1} {n \choose i} a^{n-i} b^i.
\end{array}
\end{multline*}

\item{ } In order to prove Equality~(\ref{comtet2}) we first make on the first sum the change of 
index $k = m+j$, on the second sum the change of index $j=n-k$, and we define $N = n-m$, 
so that Equality~(\ref{comtet2}) becomes
$$
\sum_{0 \leq j \leq N} {m-1+j \choose m-1} x^m (1-x)^j =
 \sum_{0 \leq j \leq N} {N+m \choose N+m-j} x^{N+m-j} (1-x)^j
$$
or, equivalently,
\begin{equation}\label{comtet3}
\sum_{0 \leq j \leq N} {m-1+j \choose m-1} (1-x)^j =
 \sum_{0 \leq j \leq N} {N+m \choose j} x^{N-j} (1-x)^j.
 \end{equation}
 
Define, for $m \geq 1$,
$$
f(m,N) = \sum_{0 \leq j \leq N} {m-1+j \choose m-1} (1-x)^j  
\ \ \mbox{\rm and} \ \
g(m,N) =  \sum_{0 \leq j \leq N} {N+m \choose j} x^{N-j} (1-x)^j.
$$
By using Pascal's rule, we write for $j \geq 1$,  
$$
{j+k \choose j} = {j+k-1 \choose j} + {j+k-1 \choose j-1} \ \mbox{\rm and} \
{N+1+j \choose j} = {N+j \choose j} + {N+j \choose j-1},
$$
so that
$$
\begin{array}{lll}
f(j+1,N) &=& (1-x) f(j+1,N-1) + f(j,N) \\
\\
g(j+1,N) &=& g(j,N) + (1-x) g(j+1,N-1).
\end{array}
$$
This implies
\begin{multline}\label{telesc}
\sum_{1 \leq j \leq m} (f(j+1,N) - f(j,N)) - \sum_{1 \leq j \leq m} (g(j+1,N) - g(j,N)) \\
= (1-x) \sum_{1 \leq j \leq m} (f(j+1,N-1) - g(j+1,N-1).
\end{multline}
But it is straightforward that
$$
\begin{array}{lll}
f(1,N) - g(1,N) &=& \displaystyle\sum_{0 \leq j \leq N} {j \choose 0} (1-x)^j 
- \sum_{0 \leq j \leq N} {N+1 \choose j} x^{N-j} (1-x)^j \\
\\
&=& \displaystyle\frac{1 - (1-x)^{N+1}}{x} - \left(\frac{1^{N+1} - (1-x)^{N+1}}{x}\right) = 0.\\
\end{array}
$$
So that Equality~(\ref{telesc}) implies
$$
f(m+1,N) - g(m+1,N) = (1-x) \sum_{1 \leq j \leq m} (f(j+1,N-1) - g(j+1,N-1)).
$$
An easy induction on $N$ then shows that, for all $N \geq 0$ and for all $m \geq 1$,
one has $f(m,N) = g(m,N)$ (the initial case $N=0$ is trivial), which proves Equality~(\ref{comtet3})
hence Equality~(\ref{comtet2}). \endpf

\end{itemize}

Two immediate corollaries are given below.

\begin{corollary}\label{2and-4}
Let $n$ be an integer $\geq 0$. Then
\begin{equation}\label{cor1-1}
\sum_{0 \leq j \leq n} 2^j {3n+1 \choose n-j} = 
(n+1) {3n+1 \choose 2n} \int_0^1 (3-2x)^n x^{2n} dx
\end{equation}
and
\begin{equation}\label{cor1-2}
\sum_{0 \leq j \leq 2n} (-4)^j {3n+1 \choose n+1+j} = \frac{1}{2} 
(n+1) {3n+1 \choose 2n} \int_{-1/2}^{3/2} (3-2x)^n x^{2n} dx.
\end{equation}
\end{corollary}

\proof First replace in Equality~(\ref{comtet1}) $n$ by $3n+1$ and $k$ by $2n$.
Then take $a = 2$ and $b = 1$. Making in the sum the change of index $i = n-j$ and
in the integral the change of variable $3-2t=2x$ yields Equality~(\ref{cor1-1}). Then take
$a = 4$ and $b = -1$. Making in the sum the change of index $j=2n-i$ and in the integral
the change of variable $t=2x$ yields Equality~(\ref{cor1-2}). \endpf

\begin{corollary}\label{eqbin}
Let $n$ be an integer $\geq 0$. Then
\begin{equation}\label{2and-3}
\sum_{0 \leq j \leq n} {3n-j \choose 2n} (1-x)^{n-j} = 
\sum_{0 \leq j \leq n} {3n+1 \choose n-j} x^j (1-x)^{n-j}.
\end{equation}
and
\begin{equation}\label{-3and-4}
\sum_{0 \leq j \leq 2n} {3n-j \choose n} (1-x)^{2n-j} =
\sum_{0 \leq k \leq 2n} {3n+1 \choose n+1+k} x^k (1-x)^{2n-k}.
\end{equation}
\end{corollary}

\proof In Equality~(\ref{comtet3}) take $m = 2N+1$. Now make in the first sum
in this Equality the change of index  $\ell = N-j$, thus obtaining
$$
\sum_{0 \leq \ell \leq N}  {3N-\ell \choose 2N} (1-x)^{N-\ell} 
= \sum_{0 \leq j \leq N} {3N+1 \choose n-\ell} x^{\ell} (1-x)^{N-\ell}.
$$
Now replace in Equality~(\ref{comtet3}) $N$ by $2N$ and take $m=N+1$, thus obtaining
$$
\sum_{0 \leq j \leq 2N} {N+j \choose N} (1-x) ^j  = 
\sum_{0 \leq j \leq 2N} {3N +1 \choose j} x^{2N-j} (1-x)^j.
$$
It remains to make in both sums the change of index $j = 2N-k$ to obtain
$$
\sum_{0 \leq k \leq 2N} {3N-k \choose N} (1-x)^{2N-k} =
\sum_{0 \leq k \leq 2N} {3N+1 \choose 2N-k} x^k (1-x)^{2N-k}
$$
which is exactly Equality~(\ref{-3and-4}) since 
$\displaystyle{3N+1 \choose 2N-k} = \displaystyle{3N+1 \choose N+1+k}$.

\begin{remark}
In \cite{alzer-prodinger} the authors define four polynomials
$$
\begin{array}{lll}
A_n(x) &= \displaystyle\sum_{0 \leq j \leq n} {3n-j \choose 2n} x^j, \ \ \ \ \
B_n(x) &= \displaystyle\sum_{0 \leq j \leq n} {3n+1 \choose n-j} x^j, \\
C_n(x) &= \displaystyle\sum_{0 \leq j \leq 2n} {3n-j \choose n} x^j, \ \ \ \ \
D_n(x) &= \displaystyle\sum_{0 \leq j \leq 2n} {3n+1 \choose n+1+j} x^j.
\end{array}
$$
After proving recurrence relations for these polynomials, the authors obtain the identities
$A_n(x+1) = B_n(x)$ and $C_n(x+1) = D_n(x)$ for which they also give direct proofs using
a variant of the Vandermonde convolution. If we compute the four quantities 
$(1+x)^{-n} A_n(x+1)$, $(1+x)^{-n} B_n(x)$, $(1+x)^{-2n} C_n(1+x)$ and $(1+x)^{-2n} D_n(x)$, 
and if we make the change of variables $x = \frac{t}{1-t}$ in our Corollary~\ref{eqbin} we obtain
$$
\begin{array}{ll}
(1+x)^{-n} A_n(x+1) &= \displaystyle\sum_{0 \leq j \leq n} {3n-j \choose 2n} (1-t)^{n-j} \\
(1+x)^{-n} B_n(x) &= \displaystyle\sum_{0 \leq j \leq n} {3n+1 \choose n-j} t^j (1-t)^{n-j} \\
(1+x)^{-2n} C_n(x+1) &= \displaystyle\sum_{0 \leq j \leq 2n} {3n-j \choose n} (1-t)^{2n-j} \\
(1+x)^{-2n} D_n(x) &= \displaystyle\sum_{0 \leq j \leq 2n} {3n+1 \choose n+1+j} t^j (1-t)^{2n-j}. \\
\end{array}
$$
Thus the relations $A_n(x+1) = B_n(x)$ and $C_n(x+1) = D_n(x)$ are exactly our
Equalities~(\ref{2and-3}) and (\ref{-3and-4}) in Corollary~\ref{eqbin}.
\end{remark}

\section{The ``mysterious'' relations given by Ruehr}

In \cite{kimura-ruehr} the following result is proved.

\begin{theorem}[Kimura-Ruehr]\label{km}

Let $f$ be a function continuous on $[-1/2, 3/2]$. Then
$$
\int_{-1/2}^{3/2} f(3x^2 - 2x^3) dx = 2 \int_0^1 f(3x^2 - 2x^3) dx
$$
\end{theorem}

Applying the result to $f(t) = t^n$, the paper announces that two combinatorial identities can be 
deduced, without giving any details. We show below that these equalities are consequences of 
Corollaries~\ref{2and-4} and \ref{eqbin} above and of Theorem~\ref{km}. Actually we even show 
that the four pairwise implied quantities are all equal (which was already proved by other authors
\cite{mtwz, alzer-prodinger}).

\begin{theorem}
We have the equalities
\begin{equation}\label{km1}
\sum_{0 \leq j \leq n} 3^j {3n-j \choose 2n}  =
\sum_{0 \leq j \leq n} 2^j {3n+1 \choose n-j} = 
\sum_{0 \leq j \leq 2n} (-4)^j {3n+1 \choose n+1+j} =
\sum_{0 \leq j \leq 2n} (-3)^j {3n-j \choose n}. 
\end{equation}
\end{theorem}

\proof
The first equality is obtained by taking $x=2/3$ in Equality~(\ref{2and-3}) of Corollary~\ref{eqbin}.
The second equality is an immediate consequence of Corollary~\ref{2and-4}. The third equality is 
obtained by taking $x = 4/3$ in Equality~(\ref{-3and-4}) of Corollary~\ref{eqbin}. \endpf

\section{Some more theory}

Relations arising in the context of distribution laws and similar to those in 
Proposition~\ref{exercises-comtet} are certainly familiar to probabilists and stasticians:  
the binomial distribution and the negative binomial distribution are respectively defined by
$$
\mbox{\rm Pr}(X = k) = {n \choose k} p^k (1-p)^{n-k} \ \ \mbox{\rm thus} \ \
\mbox{\rm Pr}(X \leq k) = \sum_{0 \leq i \leq k} {n \choose i} p^i (1-p)^{n-i} \ \
$$
and
$$
\mbox{\rm Pr}(X=k) = {k-1 \choose r-1} p^{k-r} (1-p)^r \ \ \mbox{\rm thus} \ \
\mbox{\rm Pr}(m \leq X \leq n) = \sum_{m \leq k \leq n} {k-1 \choose r-1} p^{k-r} (1-p)^r.
$$
The link with Proposition~\ref{exercises-comtet} is clear if one recalls the following two 
classical identities for the binomial and the negative binomial distributions (valid for 
$p \in [0, 1]$, $n \geq 1$, and $a \in [1,n]$),  see for example 
\cite[formulas 26.5.24, 26.5.26 on Page~945; and 6.6.1, 6.6.2 on Page~263]{abramowitz-stegun}.

\begin{proposition}[Abramowitz-Stegun]\label{a-s}
The following equalities hold
$$
\sum_{s=a}^n {n \choose s} p^s (1-p)^{n-s} = I_p(a,n-a+1)
= \frac{\int_0^p t^{a-1}(1-t)^{n-a} dt}{\int_0^1 t^{a-1}(1-t)^{n-a} dt}
$$
and
$$
\sum_{s=a}^n {n+s-1 \choose s} p^n(1-p)^s = I_{1-p}(a,n)
= \displaystyle\frac{\int_0^{1-p} t^{a-1}(1-t)^{n-1} dt}{\int_0^1 t^{a-1}(1-t)^{n-1} dt}.
$$
\end{proposition}

The reader can deduce from these equalities the statements of Proposition~\ref{exercises-comtet},
using that, if $x$ and $y$ are positive integers, then
$$
 \int_0^1 t^{x-1}(1-t)^{y-1} dt = B(x,y) =  \frac{x+y}{x} {x+y \choose x}^{-1}\cdot
 $$
 The notation $I_p(x,y)$ above stands for the {\it regularized beta function}
 $$
 I_p(x,y) = \frac{\int_0^p t^{x-1} (1-t)^{y-1} dt}{B(x,y)} = 
 \frac{\int_0^p t^{x-1} (1-t)^{y-1} dt}{\int_0^1 t^{x-1} (1-t)^{y-1} dt}
 $$
 which is the ratio of the {\it incomplete beta function} $B_p(x,y) = \int_0^p t^{x-1} (1-t)^{y-1} dt$
 to  the beta function $B(x,y) = \int_0^1 t^{x-1} (1-t)^{y-1} dt = B_1(x,y)$. To learn more about the
 incomplete beta function and its history (including formulas in Proposition~\ref{exercises-comtet}),
 the reader can consult the nice paper \cite{dutka}.
 
 \section{A link with a proof of a very weak form of the $(3x+1)$-conjecture}
 
In this section we would like to point out an application to (a generalization of) the 
$(3x+1)$-conjecture, aka.\ the Syracuse-Kakutani-Collatz-Thwaites conjecture. The 
generalization of this conjecture is the study of the function defined on the positive integers by 
 $$
 g(\ell) =
 \begin{cases}
 \displaystyle\frac{\ell}{d} \ &\mbox{\rm if $\ell \equiv 0 \bmod d$} \\
 \\
 \displaystyle\frac{n \ell - \varphi(n \ell)}{d} &\mbox{\rm otherwise} \\
 \end{cases}
 $$
where $n \geq 1$ and $d \geq 2$ are given integers such that $\gcd(d,n) = 1$, and where
$\varphi$ is the canonical surjection from the integers to a complete set ${\mathcal R}$ of 
residues modulo $d$. The conjecture (still open, even in the classical case $n=3$, $d=2$,
and ${\mathcal R} = \{0, -1\}$) asserts that, if $n < \frac{d}{d-1}$, then there exists a finite set 
${\mathcal F}$ of finite orbits under $g$ such that for each integer $\ell$, the orbit 
$\{\ell, g(\ell), g^{(2)}(\ell), \ldots\}$ of $\ell$ under $g$ is ultimately equal to an element of 
${\mathcal F}$. In other words the set $\{\ell; \ \exists k, \ g^{(k)}(\ell) \geq \ell\}$ is conjectured to be 
finite. (There is an... ultimate book by Lagarias \cite{lagarias} about the $(3x+1)$ conjecture.)
A very weak version of the conjecture has been proved, namely that the complement of this set has 
natural density zero. In one of the proofs \cite{allouche1} the author needs the following estimation 
\cite[Lemma~4, page~9-11]{allouche1}
$$
\forall \varepsilon \in (0,1), \ \exists \eta \in (0,1), \ 
\frac{1}{d^k} \sum_{|i - \frac{d-1}{d}k| > \varepsilon k} {k \choose i} (d-1)^i \leq \eta^k.
$$
In order to prove this claim, G. Tenenbaum \cite{tenenbaum} (I forgot to mention his name in 
\cite{allouche1}) suggested that I could use the equality
$$
\sum_{0 \leq i \leq m} {k \choose i} (d-1)^i = (k-m) {m \choose k} \int_{d-1}^d t^m (d-t)^{k-m-1} dt
$$
valid for $m \geq 1$ and $k \leq m-1$. The reader will have recognized (up to notation) Comtet's 
Equality~(\ref{comtet1}) of Proposition~\ref{exercises-comtet}.
 
\section{Conclusion}

The subject of sums involving binomial coefficients is huge. We have been interested here in a 
special type of sums that, although very particular, appear in many situations. Coming back
to the initial question, namely how to deduce in a simple way the Ruehr identities with binomial
sums from the Kimura-Ruehr identity for definite integrals, it would be interesting to see whether
generalizations of the Kimura-Ruehr identity (see \cite{allouche2}) could provide new identities
involving binomial sums.

\end{document}